\newif\ifprocl 
 \proclfalse 

\ifprocl
\documentclass{proc-l}
\else
\documentclass{amsart}
\fi

\usepackage{amssymb, amsmath}
\usepackage[colorlinks,citecolor=blue,urlcolor=blue]{hyperref}

\usepackage{graphicx}

\DeclareMathOperator{\dive}{div}
\DeclareMathOperator{\II}{II}
\DeclareMathOperator{\tr}{tr}
\DeclareMathOperator{\KL}{KL}
\DeclareMathOperator{\enc}{enc}
\newcommand{\dd}[1]{\mathrm{d}#1}

\newtheorem{theorem}{Theorem}[section]

\newtheorem{corollary}[theorem]{Corollary}

\theoremstyle{definition}
\newtheorem{definition}[theorem]{Definition}

\theoremstyle{remark}

\numberwithin{equation}{section}

\begin{document}

\title{Fingerprints of feedback cycles in stationary diffusions}


\ifprocl
\author{Authors}
\address{address}
\curraddr{}
\email{Email1}
\thanks{}

\else
\author[R. Schwank]{Richard Schwank$^{1,2}$}
\address{$^1$Technical University of Munich, Germany; TUM School of Computation, Information and Technology, Department of Mathematics}
\address{$^2$Konrad Zuse School of Excellence in reliable AI - relAI}
\address{$^3$Munich Center for Machine Learning}
\curraddr{}
\email{richard.schwank@tum.de}
\thanks{}

\author[M. Drton]{Mathias Drton$^{1,2,3}$}
\address{}
\curraddr{}
\email{}
\thanks{}

\fi


\date{}

\dedicatory{}

\ifprocl
\commby{Editor's name}
\fi

\begin{abstract}
When does the stationary density of a $d$-dimensional diffusion process driven by standard Brownian motion uniquely determine its drift? In general, it does not: the $d$ drift components enter the stationary Fokker--Planck equation only through a single scalar divergence. We consider drifts whose component functions only depend on a subset of their input variables, and encode this dependence with a directed graph. We exclude 2-cycles, which admit rotations that leave the stationary density invariant, but otherwise allow arbitrary graphs, including those with feedback cycles. We show that if the stationary density is sufficiently $s$-concave, e.g., strictly log-concave, then the drift is uniquely determined among graph-structured drifts in a weighted Sobolev space. For drifts with linear self-regulation, a single locally strictly convex level set of the stationary density is sufficient to identify the drift in the region it encloses. Both results follow from Raychaudhuri equations. Excluding 2-cycles allows treating nontrivial rotational and symmetric terms simultaneously, extending the scope of the corresponding Bochner- and Reilly-type arguments. We illustrate the results by recovering the Hill-type regulation of a stochastic genetic repressilator from samples of its stationary distribution.
\end{abstract}

\ifprocl
\else
\fi

\maketitle

\section{Introduction}
In many settings of practical interest, a stochastic process is not observed continuously. Observation frequencies may be physically limited, or measurement may destroy the sample. From an engineering perspective, only the final process state may be of interest, like in generative models. For these settings, algorithms have been designed that learn a diffusion process which approximately has a given data distribution as its stationary distribution \cite{varando20a, Song2021GenerativeSDE, BotvinickGreenhouse2023, Dettling2024Lasso, Lorch2024StationaryDiffusions}.  

A relevant theoretical question is in what sense there is a \emph{unique} diffusion process \begin{equation}\label{eq:sde_const_diff}
    \dd{x}(t) = b(x(t))\,\dd{t} + \dd{W}(t)
\end{equation} with stationary density $\rho$. Here, $b$ denotes the drift vector field, $W$ is a $d$-dimensional Wiener process, and the equation is interpreted in Itô sense. Under mild regularity conditions, a stationary density $\rho$ of \eqref{eq:sde_const_diff} solves the stationary Fokker--Planck equation \cite{Bogachev2015FPBook}\begin{equation}\label{eq:FP_const_diff}
    \dive(b\rho) = \tfrac{1}{2}\Delta\rho.
\end{equation} This is a single scalar equation for the $d$ component functions of $b$. Commonly, for example in generative machine learning, attention is restricted to reversible processes, i.e., to gradient fields $b$ \cite{bleile2026efficient}. 
In this case the stationary density $\rho$ essentially determines the drift $b$ uniquely, for example $b=\tfrac{1}{2}\nabla\log \rho$ in \eqref{eq:sde_const_diff}. For non-reversible processes, however, the drift is in general not unique: if $v$ is any vector field with $\dive(v\rho)=0$, then $b+v$ solves \eqref{eq:FP_const_diff} as well.

We consider drifts $b$ with low-dimensional structure. Precisely, we assume that the component functions $b^i(x)$ only depend on input variables $(x_j)_{j\in \mathrm{pa}(i)}$ for some subsets $\mathrm{pa}(i)\subseteq \{1,\ldots, d\}$. We depict these dependencies by a directed graph $\mathcal{G}$, in which node $i$ has parents $\mathrm{pa}(i)$. Such systems are found in probability theory \cite{Lacker2023StationaryDiffusionOnTree}, computational biology \cite{aguirre2025gene},  statistics \cite{Gregorio2025DiffusionProcessRegression},   and causality \cite{boeken2026causalgraphsmarkovproperties}. Usually, the associated stochastic processes are non-reversible. 

Previously, uniqueness had only been shown without feedback cycles \cite{schwank2026arxiv} or among linear drifts, i.e., for Ornstein--Uhlenbeck processes \cite{Dettling2023}. The difficulty with cyclic graphs $\mathcal{G}$ is that certain rotations leave the stationary distribution invariant. For example,  
the vector field $v = \rho^{-1}(\partial_i\rho\, e_j - \partial_j\rho\, e_i)$ satisfies $\dive(v\rho)=0$ for any distinct $i,j$. If $\rho$ has independent coordinates, then $v_i$ depends only on $x_j$ and $v_j$ only on $x_i$, so $b+v$ is again $\mathcal{G}$-structured whenever $\mathcal{G}$ contains the $2$-cycle $i\leftrightarrows j$. These rotations are observed in practice \cite[Fig.~3]{BotvinickGreenhouse2023}. We show that for a large class of densities $\rho$, any undetectable rotation in fact requires a $2$-cycle.

Concretely, we prove that the stationary density $\rho$ uniquely determines the drift when the structure graph $\mathcal{G}$ has no $2$-cycles and when $\rho$ satisfies a suitable convexity condition. A special case of our later results reads as follows: \begin{theorem}\label{thm:illustrate}
    Let $\mathcal{G}$ be a directed graph on $d\geq 1$ nodes without $2$-cycles and let $\rho\in C^2(\mathbb{R}^d)$ be a positive probability density. Assume $\rho$ is  \emph{strictly log-concave} in the sense that $\nabla^2\log\rho$ is negative definite. Consider the drift class \begin{equation*}
        \mathcal{B}(\mathcal{G}) := \{b\in C^2(\mathbb{R}^d,\mathbb{R}^d) \mid b \text{ is }\mathcal{G}\text{-structured and }b, \nabla b\in L^2(\rho)\}.
    \end{equation*} Then, there is at most one drift $b\in\mathcal{B}(\mathcal{G})$ such that $\rho$ is a stationary density of the process \eqref{eq:sde_const_diff} in the sense of \eqref{eq:FP_const_diff}. 
\end{theorem} 
The main results in Theorem \ref{thm:global} and \ref{thm:local_identi} extend Theorem \ref{thm:illustrate} in two ways. First, we weaken log-concavity to $s$-concavity and we extend the constant diffusivity from \eqref{eq:sde_const_diff} to arbitrary, known diffusivity. Second, we prove local uniqueness on a subset of $\mathbb{R}^d$ under a partial linearity assumption and  a local level-set convexity assumption. We numerically illustrate our results in Figure \ref{fig:repressilator} for the repressilator \cite{Elowitz2000}, which is a well-known model for cyclic gene regulation.

Mathematically, we frame identifiability as an inverse problem for the stationary Fokker--Planck equation and link uniqueness with non-existence of certain weighted divergence-free vector fields. We show global drift identifiability by a Bochner-type argument \cite{Bochner1946} featuring the Bakry--\'Emery curvature induced by the density $\rho$, and we show local identifiability by a Reilly-type argument \cite{Raulot2011}. These results rely on the surprising link that absence of $2$-cycles supplies exactly the sign condition needed to apply both arguments, which usually require restrictions to smaller classes such as gradient fields. We motivate and present the proof from the Raychaudhuri equation and its weighted version, which also supplies a physical interpretation.   

Our results are of interest in several ways. We find it remarkable that a single first-order divergence equation can uniquely determine a $d$-dimensional vector field among a large class of nonlinear alternatives. Moreover, the graph structure is a non-standard sufficient condition for the geometric arguments we use. Our results also provide identifiability guarantees for methods that estimate diffusion processes from stationary data \cite{BotvinickGreenhouse2023, Lorch2024StationaryDiffusions}. For experimental design, our theory suggests measuring mediators: if $j$ regulates $i$ only through a third variable $k$, then observing $k$ turns the $2$-cycle $i\leftrightarrows j$ into the $3$-cycle $i\to j\to k\to i$, to which our results apply.

The remainder of the paper derives the Raychaudhuri equations in Section \ref{sec:backgr}, presents the main results in Section \ref{sec:theory}, illustrates the theory numerically in Section \ref{sec:numerics}, and discusses further related work in Section \ref{sec:rel_work}. 

\subsection*{Notation}
We write the Euclidean coordinates of vector and matrix fields as $w_i$ and $A_{ij}$, respectively.  
We write $\langle\cdot,\cdot\rangle$ for the Euclidean product and $|w|$, $|A|$ for the induced vector and Frobenius norm.
We write $\nabla_i = \partial/\partial x_i$, so that $\dive w = \sum_i\nabla_i w_i$ is the divergence, $\nabla_w = \sum_k w_k\nabla_k$ the directional derivative along $w$, and $(\nabla w)_{ij} = \nabla_i w_j$ the Jacobian of a vector field $w$. For a scalar function $\gamma$, $\nabla \gamma$ denotes the gradient, $\nabla^2 \gamma$ the Hessian, and $\Delta \gamma = \dive \nabla \gamma$ the Laplacian.

\section{Background: The Raychaudhuri equation}\label{sec:backgr}
The purpose of this section is to compactly derive the Raychaudhuri equation \eqref{eq:raychaudhuri} and a probability-weighted version \eqref{eq:Raychaudhuri_wghtd} in Euclidean space. The second version features a weighted divergence \eqref{eq:weighted_div_def}, which introduces curvature-like terms.

Let $w$ be a $C^2$-vector field in $d$-dimensional Euclidean space. The Raychaudhuri equation concerns $\nabla_w \dive w$, meaning how the divergence $\dive w$ changes along the flow generated by $w$ itself. In physics, more concretely in general relativity, the equation describes local expansion along flows in spacetime \cite{Kar2007}. 
To derive the equation, first consider the following counterpart: \begin{equation*}
    \dive (\nabla_w w) = \sum_{i,k\leq d}\nabla_i(w_k \nabla_k w_i) = \tr((\nabla w)^2) + \sum_{i,k\leq d}w_k\nabla_i\nabla_k w_i,
\end{equation*} where $\tr$ denotes the matrix trace. On Euclidean space, we may exchange  $\nabla_i$ and $\nabla_k$, and the final term is exactly $\nabla_w \dive w$. 

We have obtained the \emph{Raychaudhuri equation}: \begin{equation}\label{eq:raychaudhuri}
    \nabla_w \dive w = \dive(\nabla_w w) - \tr((\nabla w)^2). 
\end{equation}

Our main theorem requires a version of this equation under an additional source of curvature. Note that exchanging derivatives on a curved space would have led to an additive nonzero curvature term governed by the traced Ricci identity. 
The curvature we consider does not originate from the space itself, but from a probability density $\rho$ on $\mathbb{R}^d$, which we assume to be $C^2$ and  positive everywhere. Setting $\varphi := -\log \rho$, we define a \emph{weighted divergence}: \begin{equation}\label{eq:weighted_div_def}
    \dive_\varphi w := \dive(w) - \langle w, \nabla\varphi\rangle = \frac{\dive(w\rho)}{\rho}.
\end{equation} Note that the numerator also appears in the stationary Fokker--Planck equation \eqref{eq:FP}. One motivation for this exact definition is that $\dive_\varphi$ satisfies a clean divergence theorem when integrating over domains weighted by $\rho$. 

To obtain a Raychaudhuri equation for the weighted divergence, note that the left hand side of equation \eqref{eq:raychaudhuri} misses $\nabla_w\langle w,\nabla\varphi\rangle$, which decomposes into $\langle\nabla_w w, \nabla\varphi\rangle$ and $\nabla^2\varphi(w,w) :=\langle w, (\nabla^2 \varphi) w\rangle $. The first term is what completes the weighted divergence term on the right hand side. The Hessian term is commonly absorbed into any Ricci terms, resulting in the \emph{Bakry--\'Emery Ricci tensor} \cite{BGL2014Book}. To summarize, we have the desired \emph{weighted Raychaudhuri equation}: \begin{equation}\label{eq:Raychaudhuri_wghtd}
    \nabla_w \dive_\varphi w = \dive_\varphi(\nabla_w w) - \tr((\nabla w)^2) - \nabla^2\varphi(w,w).
\end{equation}

\section{Theoretical results}\label{sec:theory}
We observe a probability density $\rho\in C^2(\mathbb{R}^d)$, which we assume to be positive everywhere, and again set $\varphi:=-\log \rho$.
If $\rho$ is the stationary or ergodic density of a diffusion process, it must satisfy the stationary Fokker--Planck equation under mild regularity conditions. We therefore frame  the inverse problem as follows.

\begin{definition}\label{def:identif}
    Let $\mathcal{B}$ be a class of $C^2$ vector fields on $\mathbb{R}^d$.  We say $\rho$ \emph{identifies the drift relative to the diffusivity among} $\mathcal{B}$\emph{-vector fields}, if for all $C^2$ matrix fields $\Sigma$ 
    there is at most one $b\in\mathcal{B}$ satisfying the stationary Fokker--Planck equation \begin{equation}\label{eq:FP}
    \dive(b \rho) = \frac{1}{2}\sum_{i,j\leq d}\nabla_i \nabla_j ( (\Sigma\Sigma^T)_{ij}\, \rho). 
\end{equation} 
\end{definition}

If $\mathcal{B}$ is closed under subtraction, identifiability is equivalent to the fact that \begin{equation}\label{eq:identif_equiv_cond}
    \dive_\varphi v = 0 \text{ implies } v=0 \text{ for all }v\in\mathcal{B},
\end{equation} where $\dive_\varphi v = \rho^{-1}\dive(v \rho)$ denotes the weighted divergence from \eqref{eq:weighted_div_def} and  where $v$ represents a difference of drifts. To show that Definition \ref{def:identif} implies \eqref{eq:identif_equiv_cond} by contraposition, we may use $\Sigma\equiv 0$ or $\Sigma = \rho^{-1/2}I$ with the identity matrix $I$. The vector field class $\mathcal{B}$ we consider builds on the following graph terminology.

\begin{definition}
    A directed graph $\mathcal{G}$ \emph{has no} $2$\emph{-cycles} if for all distinct nodes $i,j$ the graph $\mathcal{G}$ contains at most one of the edges $i\to j$ or $j\to i$. Self-loops $i\to i$ are allowed. The $C^2$ vector field $b$ is $\mathcal{G}$\emph{-structured} according to a directed graph $\mathcal{G}$ if $\nabla_i b_j = 0$  whenever the edge $i\to j$ is not contained in $\mathcal{G}$. 
\end{definition}

Fix a graph $\mathcal{G}$ without $2$-cycles and consider the class of $\mathcal{G}$-structured vector fields. The class is closed under subtraction, hence the divergence characterization \eqref{eq:identif_equiv_cond} applies.

One motivation for studying this class is the following lower bound used in \cite{Dettling2023}: When $v$ has a structure graph $\mathcal{G}$ without $2$-cycles, then \begin{equation}\label{eq:trace_nonneg}
    \tr((\nabla v)^2) = \sum_{i,j\leq d}\nabla_i v_j \nabla_j v_i = \sum_{i\leq d} (\nabla_{i} v_i)^2 \geq 0,
\end{equation} since at least one of $\nabla_i v_j$ or $\nabla_j v_i$ is zero for any $i\neq j$ due to absence of $2$-cycles.

This bound has a concrete interpretation in the context of the Raychaudhuri equation. Decompose $\nabla v = K + S$, where $K_{ij} := \tfrac{1}{2}(\nabla_i v_j - \nabla_j v_i)$ is the anti-symmetric part 
corresponding to rotation in the flow generated by $v$, and where $S_{ij} := \tfrac{1}{2}(\nabla_i v_j + \nabla_j v_i)$ is the symmetric part 
combining shear and expansion. 
Then, \begin{equation*}
    \tr((\nabla v)^2) = \sum_{i,j\leq d}\nabla_i v_j \nabla_j v_i = \sum_{i,j\leq d}\left(S_{ij}S_{ij} - K_{ij} K_{ij}\right) = |S|^2 - |K|^2.
\end{equation*} Hence the identity \eqref{eq:trace_nonneg} means that the symmetric part of $\nabla v$ dominates the antisymmetric part in terms of magnitude for any $v$ with $2$-cycle-free structure graph. Note that $\dive v$ only depends on the symmetric part of $\nabla v$ 
and that the divergence is the only term where derivatives of $v$ enter the condition \eqref{eq:identif_equiv_cond}. This illustrates why a control of $K$ via $S$ may be useful.

\subsection{Global identifiability}
To apply the criterion \eqref{eq:identif_equiv_cond}, consider a vector field $v$ with $\dive_\varphi v =0$. The weighted Raychaudhuri equation \eqref{eq:Raychaudhuri_wghtd} then reads \begin{equation}\label{eq:wghtd_raych_red}
    \dive_\varphi(\nabla_v v) = \tr((\nabla v)^2) + \nabla^2\varphi(v,v).
\end{equation} Our proof strategy relies on making the right hand signed. A starting point is to use nonnegativity of the trace term  \eqref{eq:trace_nonneg} in combination with positive definiteness of $\nabla^2 \varphi$, i.e., assuming strict log-concavity of $\rho$ as in Theorem \ref{thm:illustrate}. However, this approach can be strengthened by reconsidering  the final step of \eqref{eq:trace_nonneg}. For any vector field $v$ with $2$-cycle-free structure graph and $\dive_\varphi v = 0$, convexity of the square and $\dive  v = \langle v, \nabla\varphi\rangle$ imply that \begin{equation}\label{eq:s-conc_lower_bd}
    \tr((\nabla v)^2) + \nabla^2\varphi(v,v) \geq \tfrac{1}{d} \langle v, \nabla \varphi\rangle^2 + \nabla^2\varphi(v,v) = d\, \rho^{1/d}\, \nabla^2 (\rho^{-1/d})(v,v).
\end{equation}
The middle term is also known as the $N=0$ dimensional Bakry--\'Emery curvature \cite{Kolesnikov2018}, 
while the final Hessian demonstrates the connection with $s$-concave densities \cite{Saumard2014LogConcaveReview}. 
\begin{definition}
    Let $s\leq 0$ and $M\subset\mathbb{R}^d$ be convex. The density $\rho$ is \emph{strictly} $s$\emph{-concave on} $M$ if $\nabla^2 (\rho^s)\succ 0$ on $M$ for $s<0$ and if  $\nabla^2\varphi\succ0$ for $s=0$. 
\end{definition}
In the definition, $A \succ 0$ denotes positive definiteness. Note that if $\rho$ is strictly $s$-concave, then it is strictly $s'$-concave for all $s' < s\leq 0$. For example, the multivariate $t$-density proportional to \begin{equation*}
    (1 + \langle x, C x\rangle)^{-\tfrac{\nu + d}{2}}
\end{equation*} for some $\nu >0, C\succ0$ is strictly $-(\nu + d)^{-1}$-concave, 
in particular $-d^{-1}$-concave, 
but not log-concave. We can now prove the main result.

\begin{theorem}\label{thm:global}
    Let $\mathcal{G}$ be a graph without $2$-cycles on $d$ nodes. Consider the class
    \begin{equation*}
        \mathcal{B}(\mathcal{G}) := \{b \mid b \text{ is }\mathcal{G}\text{-structured and }b, \nabla b\in L^2(\rho)\}.
    \end{equation*} When the density $\rho$ is strictly $-d^{-1}$-concave on $\mathbb{R}^d$, then $\rho$ identifies the drift relative to the diffusivity among $\mathcal{B}(\mathcal{G})$-vector fields.
\end{theorem}
\begin{proof}
    Since $\mathcal{B}(\mathcal{G})$ is closed under subtraction, we target the criterion \eqref{eq:identif_equiv_cond}. Let $v\in\mathcal{B}(\mathcal{G})$ with $\dive_\varphi v = 0$, let $r>0$, and set $B_r := \{x\in\mathbb{R}^d\ : |x|\leq r\}$. Integrating \eqref{eq:wghtd_raych_red} with respect to $\rho$ yields that \begin{equation*}
        \int_{B_r} \tr((\nabla v)^2) + \nabla^2\varphi(v,v) \,\dd{\rho} = \int_{B_r} \dive_\varphi (\nabla_v v)\,\dd{\rho} = \int_{\partial B_r} \langle \nabla_v v, n\rangle \rho \,\dd{S}=: e_r.
    \end{equation*} Using the co-area formula and $|v|, |\nabla v|\in L^2(\rho)$, it follows that \begin{equation*}
        \int_0^\infty |e_r|\,\dd{r} \leq \int_{\mathbb{R}^d} \frac{1}{2}(|v|^2 + |\nabla v|^2)\,\dd{\rho} < \infty,
    \end{equation*} where we applied both Cauchy-Schwarz and Young's inequality. This result implies $\liminf_{r\to\infty} |e_r| = 0$. Using $\dive_\varphi v =0$, absence of $2$-cycles, $-d^{-1}$-concavity of $\rho$, and monotone convergence, the argument \eqref{eq:s-conc_lower_bd} shows that \begin{equation*}
        0 \leq \int_{\mathbb{R}^d} d\, \rho^{1/d}\, \nabla^2 (\rho^{-1/d})(v,v)\,\dd{\rho} \leq  \int_{\mathbb{R}^d} \tr((\nabla v)^2) + \nabla^2\varphi(v,v) \,\dd{\rho} = 0.
    \end{equation*} Since $\rho >0$, it follows that $\nabla^2 (\rho^{-1/d})(v,v) =0$ everywhere by continuity. 
    Strict  $-d^{-1}$-concavity of $\rho$ implies $v=0$ as required by \eqref{eq:identif_equiv_cond}. 
\end{proof}
This theorem directly generalizes Theorem 32  in \cite{Dettling2023}, which treats Gaussian distributions and proves that the stationary density identifies the drift relative to constant diffusivities among \emph{linear} vector fields. Our Theorem applied to Gaussian distributions also covers nonlinear drift alternatives and non-constant diffusivities. 

Theorem \ref{thm:global} is mostly relevant in cases where $\mathcal{G}$ has a cycle. Otherwise, identifiability follows already without the convexity assumption on $\rho$ \cite{schwank2026arxiv}.

The condition that the drift $b$ and its Jacobian $\nabla b$ are square integrable with respect to the stationary density $\rho$ is mild, especially for light-tailed $\rho$. Such integrability conditions are common in the literature to establish regularity results and prevent explosions \cite{Haesung2021FullSDEPackageUnderLocalIntegrability, schwank2026arxiv}. 

Note that the threshold $d^{-1}$ can be improved to the inverse number of self-loops in $\mathcal{G}$. This has an interpretation in terms of finite-dimensional Bakry--\'Emery curvature, but we did not include the improvement for a lack of statistical examples.

\subsection{Linear self-regulation} In this section, we present a local identifiability result from a local convexity condition on the density $\rho$.

\begin{definition}\label{def:identif_loc}
    Let $\mathcal{B}$ be a class of $C^2$ vector fields and $S\subset \mathbb{R}^d$. We say $\rho$ \emph{identifies the drift relative to the diffusivity among $\mathcal{B}$-vector fields locally on $S$}, if for any $C^2$ matrix field $\Sigma$, all drifts $b\in\mathcal{B}$ solving the Fokker--Planck equation \eqref{eq:FP} are identical on $S$. 
\end{definition}
If $\mathcal{B}$ is closed under subtraction, Definition \ref{def:identif_loc} is equivalent to the fact that \begin{equation}\label{eq:identif_equiv_cond_loc}
    \dive_\varphi v = 0 \text{ implies } v=0 \text{ on } S \text{ for all }v\in\mathcal{B}.
\end{equation} Let $\mathcal{G}$ be a graph without $2$-cycles. In this section, we consider the following class:
\begin{equation}\label{eq:def_B_lin}
    \mathcal{B}_\text{lin}(\mathcal{G}) := \{b \mid b \text{ is } \mathcal{G}\text{-structured and } \exists\, \alpha\in\mathbb{R}^d : \nabla_{i} b_{i} = \alpha_{i} \text{ for all } i\leq d \},
\end{equation} Integrating the defining condition shows that the $i$-th component $b_i$ is of the form $\alpha_i x_i + r_i$ for some function $r_i$ only depending on the proper parents of node $i$ in $\mathcal{G}$. In other words, the drifts in the class $\mathcal{B}_\text{lin}(\mathcal{G})$ have linear self-regulation, which is a common feature of diffusion models. 

Since $\mathcal{B}_\text{lin}(\mathcal{G})$ is closed under subtraction, consider $v\in \mathcal{B}_\text{lin}(\mathcal{G})$ with $\dive_\varphi v = 0$. Due to the linear self-regulation, it holds that \begin{equation*}
    \dive v(x) = \text{const} = \langle v(x),\nabla\varphi(x)\rangle \qquad\forall x\in\mathbb{R}^d. 
\end{equation*} If $\rho$ has a critical point $x_0$, then $\nabla\varphi(x_0)=0$ and this constant is zero. To summarise, the following two identities hold: \begin{align}\label{eq:tangent_and_zeroDive}
    \langle v, \nabla \varphi\rangle & = 0,  & \dive v & = 0.
\end{align} The left identity means that $v$ is tangential to level sets of $\varphi$, which indicates that integrating over sublevelsets of $\varphi$ might produce manageable boundary terms.  

The right identity suggests to use the unweighted Raychaudhuri equation analogously to the weighted version earlier. The advantage is that the unweighted version does not contain the Hessian of $\varphi$, reducing curvature assumptions.

Before stating the local identifiability result, recall that the second fundamental form $\II$ of a hypersurface $S$ with $C^1$ unit outward normal $n$ is defined as $\II(v,w):= \langle n, \nabla_v w\rangle$ for tangential vector fields $v,w$. We call $\II$ negative definite if $\II(v,v) < 0$ on $S$ for all tangential $v\neq 0$. In this case, $S$ locally is the graph of a strictly convex function. If $S$ additionally is connected and compact, 
there is a unique bounded connected component of $\mathbb{R}^d\setminus S$, which we denote as $\enc(S)$ \cite{Lima1988}.

\begin{theorem}\label{thm:local_identi}
    Let $\mathcal{G}$ be a graph without $2$-cycles on $d\geq 2$ nodes. Assume there is  $t > 0$ and a bounded, connected component $S$ of $\{\rho = t\}$ which is an orientable hypersurface with negative definite second fundamental form $\II$ and $\nabla\rho \neq 0$ on $S$.
    
    Then, $\rho$ identifies the drift relative to the diffusivity among $\mathcal{B}_\text{lin}(\mathcal{G})$-vector fields locally on $\enc(S)\cup S$.
\end{theorem}
\begin{proof}
    Aiming first to apply the criterion \eqref{eq:identif_equiv_cond_loc} to the surface $S$, we consider $v\in \mathcal{B}_\text{lin}(\mathcal{G})$ with $\dive_\varphi v = 0$. Abbreviate $M:=\enc(S)$ and note that 
    the outer unit normal field $n$ to $S$ is parallel to $\nabla\varphi$ since $S\subset \{\rho = t\}$ and $\nabla\rho\neq 0$.

    Further, there exists $x_0\in M$ with $\nabla\rho(x_0) = 0$: Either $\rho \equiv t$ on the compactum $\bar{M}$ or there is an extremum inside of $M$. 
    Consequently,  $v$ satisfies equation \eqref{eq:tangent_and_zeroDive} as discussed. In particular, $v$ is tangential to $S=\partial M$. Inserting $\dive v = 0$ into the unweighted Raychaudhuri equation \eqref{eq:raychaudhuri} yields that \begin{equation*}
        \dive (\nabla_v v) = \tr((\nabla v)^2).
    \end{equation*} Integrating over $M$ and applying the divergence theorem shows that \begin{equation*}
        \int_{\partial M} \II(v,v) \dd{S} = \int_{\partial M} \langle n, \nabla_v v \rangle \dd{S} = \int_M \tr((\nabla v)^2) \dd{x} = \int_M \sum_{i\leq d} \alpha_i^2 \dd{x}\geq 0,
    \end{equation*} where the decomposition of the trace is due to absence of $2$-cycles \eqref{eq:trace_nonneg} and the constants $\alpha_i$ come from linear self regulation \eqref{eq:def_B_lin}. 
    By assumption, $\II(v,v)\leq 0$ on $\partial M$, which, together with the integral lower bound, forces $\II(v,v)=0$ on $\partial M$ by continuity and $\alpha_i=0$ for all $i\leq d$. 
    
    Negative definiteness of $\II$ implies $v=0$ on $S=\partial M$. Further, for any $i\leq d$ and $x\in\enc(S)$ let $x^+$ be an intersection of the line $x + \mathbb{R}e_i$ with $S$. Since $\alpha_i=0$, the function $v_i(x)$ does not depend on $x_i$, which implies $v_i(x)= v_i(x^+) = 0$. 
\end{proof}

Remarkably, this Theorem only requires a level set piece of a certain geometry and no condition on the density $\rho$ inside the level set. Consider for example a finite Gaussian mixture in dimension $d\geq 2$ with one covariance matrix $\Sigma_{\max}$ strictly dominating all others. Since $\nabla^2\varphi \to \Sigma_{\max}^{-1}\succ 0$ as $|x|\to\infty$, the level sets $\{\rho = t\}$ satisfy the condition in Theorem \ref{thm:local_identi} as $t\searrow 0$, while the interior density landscape can be nonconvex. Letting $t\searrow 0$ thus even results in a global identifiability result for these mixture densities $\rho$. 

Another case where the local identifiability result has global consequences is with analyticity. Since $S$ encloses a nonempty set, the density identifies the drift relative to the diffusivity among $\mathcal{B}_\text{lin}(\mathcal{G})\cap \{\text{real-analytic vector fields}\}$.

Finally, an important scenario in which Theorem \ref{thm:local_identi} applies is local concavity.

\begin{corollary}\label{cor:local_max}
    Let $\mathcal{G}$ be a graph without $2$-cycles on $d\geq 2$ nodes. Assume $\rho$ is strictly $s$-concave for some $s\leq 0$ on a bounded, nonempty, open, and convex set $M\subset\mathbb{R}^d$. Assume $\rho$ attains a local maximum at $x_0\in M$. Then, $\rho$ identifies the drift relative to the diffusivity among $\mathcal{B}_\text{lin}(\mathcal{G})$-vector fields locally on $M \cap \{\rho > \max_{\partial M} \rho\}$. 
\end{corollary}
\begin{proof}

    Let $\rho(x_0) > t > \max_{\partial M} \rho$. Set $S := M \cap \{\rho = t\}$ and define $q:= \varphi$ if $s=0$ and $q:= \rho^s$ otherwise. 
    
    Note that $\nabla \rho(x) = 0 $ if and only if $ \nabla q(x) = 0$. Hence, another root of $\nabla\rho$ on $M$ beyond $x_0$ would contradict strict convexity of $q$. 
    Therefore, the implicit function theorem provides local charts of $\{\rho = t\}$. The lower bound on $t$ ensures these are also charts of $S$ and that $S = \partial(M\cap \{\rho > t\})$ is connected, as it is the boundary of a convex body in dimension $d\geq 2$. 
    An outer normal field is given by  $\nabla\varphi / |\nabla \varphi|$. Thus $S$ is a bounded, connected, orientable hypersurface contained in $\{\rho = t\}$.

    Let $w$ be a tangential field to $S$, i.e. satisfying $\langle w, \nabla \varphi\rangle = 0$. Taking $\nabla_w$ on both sides yields that \begin{equation*}
        |\nabla \varphi|\cdot \II(w,w) = \langle \nabla\varphi, \nabla_w w\rangle = - \langle w, (\nabla ^2 \varphi) w\rangle = - f(s)\langle w, (\nabla ^2 q) w\rangle,
    \end{equation*} where the final identity uses $\nabla^2\varphi = \nabla^2 q$ up to $f(s)=|s|^{-1} \rho^{|s|}$ on the orthogonal complement of $\nabla\varphi$ for $s<0$ as seen in equation \eqref{eq:s-conc_lower_bd}. Trivially, $f(0)=1$. Recalling $\nabla\varphi \neq 0$ and $\nabla^2 q \succ 0$ on $S$ shows that $\II$ is negative definite. Applying Theorem \ref{thm:local_identi} yields local identifiability on $S\cup \enc(S)$. As $t$ was arbitrary, the claim follows. 
\end{proof}

\section{Numerical illustration}\label{sec:numerics}

\begin{figure}
    \centering
    \includegraphics[width=0.45\linewidth]{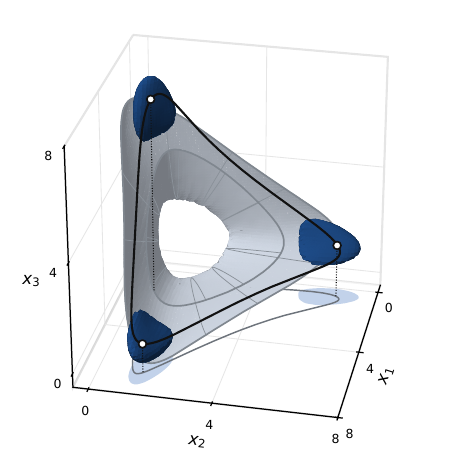}
    \raisebox{0.45cm}{
    \includegraphics[width=0.45\linewidth]{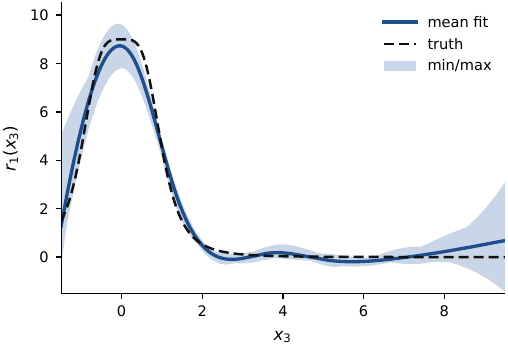}
    }
    \caption{\emph{Left}: Illustration of the repressilator's stationary density for $\alpha=-1, \beta=9,\gamma=4$, and $\sigma=0.7$. White points mark the maxima of the density, and the black reference line shows the periodic orbit of the associated deterministic system ($\sigma=0$). The blue surfaces are the level sets at 75\% of the maximal density. The grey surface is a torus-shaped level set at 34\% of the maximal density value, which is cut in half by a plane orthogonal to $(1,1,1)$. \emph{Right}: Mean fitted regulation function $r_1(x_3)$ from $10$ experiments with $N=1000$ samples of $\rho$ each, together with the pointwise minimal and maximal fit over all experiments. The fit agrees well with the true regulation function $9/(1+x_3^4)$.}
    \label{fig:repressilator}
\end{figure}

The \emph{repressilator} \cite{Elowitz2000} is a gene regulation model that produces oscillations, i.e., a ``biological clock". It demonstrates the rich dynamical behavior that feedback cycles can produce in combination with non-linearity. Here, we consider a stochastic version in three variables \cite{Shen2025}, with indices understood modulo three: \begin{equation}\label{eq:repressilator}
    \dd{x_i(t)} = \left(\alpha x_i(t) + \frac{\beta}{1 + x_{i-1}^{\gamma}(t)}\right)\dd{t} + \sigma\, \dd{W_i(t)}, \qquad i=1,2,3.
\end{equation} The system approaches a unique stationary density $\rho$ when $\gamma$ is an even integer and $-\alpha,\sigma >0$. Lipschitz continuity of the drift guarantees a global solution, the Lyapunov function $|x|^2$ ensures recurrence, 
and $\sigma>0$ provides irreducibility. 
Hence, a stationary density $\rho$ exists, 
is unique, 
and the process converges weakly to $\rho$. 
Figure \ref{fig:repressilator}
illustrates a Monte Carlo approximation of $\rho$ for $\alpha=-1, \beta=9,\gamma=4$, and $\sigma=0.7$. The density broadly concentrates around the limit cycle of the deterministic system with $\sigma=0$ \cite{MalletParet1990}. The level sets of $\rho$ are ellipsoidal near the maxima and become torus-shaped in lower density regions.

The remainder of this section treats learning the repressilator drift function from the stationary density $\rho$ from the left of Figure \ref{fig:repressilator}, both theoretically and empirically. Notice that the drift is $\mathcal{G}$-structured for the $3$-cycle $\mathcal{G}$ with self-loops, which has no $2$-cycles. Since the self-regulation term is linear, the drift lies in $\mathcal{B}_\text{lin}(\mathcal{G})$.

Local identifiability is guaranteed near the maxima of $\rho$. More precisely, numerical checks show that $\nabla^2\varphi$ is positive definite at the maxima. Hence, Corollary \ref{cor:local_max} applies and $\rho$ locally identifies the repressilator drift relative to the diffusivity among $\mathcal{B}_\text{lin}(\mathcal{G})$ drifts. Concretely, \begin{equation}\label{eq:B_lin_3cycle}
    \mathcal{B}_\text{lin}(\mathcal{G}) = \{b \mid \,  b_i(x) = \alpha_i x_i + r_i(x_{i-1})\text{ with }\alpha_i\in\mathbb{R}, r_i\in C^2(\mathbb{R})\text{ for all }i\leq 3\},
\end{equation} and any two drift vector fields $b_1, b_2\in\mathcal{B}_\text{lin}(\mathcal{G})$ for which $\dd{x(t)}=b_j(x(t))\dd{t} + 0.7 \dd{W(t)}$ has $\rho$ as a stationary distribution must agree perfectly in some neighborhood of the points marked white in Figure \ref{fig:repressilator}. In particular, 
there is exactly one real analytic drift in $\mathcal{B}_\text{lin}(\mathcal{G})$ with $\rho$ as the induced stationary distribution for $\sigma=0.7$.

Indeed, we can numerically recover the repressilator drift reasonably well from $N=1000$ approximately independent Monte Carlo samples of $\rho$ as shown on the right of Figure \ref{fig:repressilator}. We adapt the method by \cite{BotvinickGreenhouse2023}, aiming to find $b\in\mathcal{B}_\text{lin}(\mathcal{G})$ with minimal discrepancy between the data distribution and the solution $\rho$ of the stationary Fokker--Planck equation \begin{equation}\label{eq:FP_scalar_diff}
    \tfrac{\sigma^2}{2} \Delta \rho = \dive(b \rho).
\end{equation} We employ cubic splines to parametrize the regulation functions $r_i$ introduced in \eqref{eq:B_lin_3cycle}; further implementation details are discussed below. Figure \ref{fig:repressilator} shows that the average estimated $r_1$ over 10 independent repetitions accurately recovers the ground truth, which is consistent with our theoretical result. Note that global recovery was possible empirically even among a non-analytic function class in this example. 

In more detail, the numerical method first discretizes a domain containing the training data into cubes of side length $h=0.25$, and then computes a histogram $\hat{\rho}_h$. We parametrize $\mathcal{B}_\text{lin}(\mathcal{G})$ through $(\alpha_i, r_i)_{i=1,2,3}$, and each function $r_i$ by a cubic spline with $10$ basis functions, 
giving 33 parameters $\theta$ in total. Write $\rho_h(\theta)$ for the finite volume Scharfetter-Gummel solution of \eqref{eq:FP_scalar_diff} with drift $b(\theta)$ and $\sigma=0.7$.  
The strategy is to optimize the KL-divergence $\KL(\hat{\rho}_h \| \rho_h(\theta))$ over $\theta$ with the quasi-Newton method L-BFGS-B with box constraints. This requires $\nabla_\theta \rho_h(\theta)$, which when computed naively would require solving one linear equation on the scale of the finite volume system per parameter. The method instead exploits that $\nabla_\theta \rho_h(\theta)$ is only needed in the pairing $(\partial_{\rho_h}\KL)^T\nabla_\theta \rho_h(\theta)$, for which solving one adjoint problem suffices. 
To avoid instabilities with the KL-divergence when $\hat{\rho}_h$ has mass in cells where $\rho_h$ is almost zero, we perform a warm-up using the $L^2$-loss before optimizing the KL-loss. Code and figures for this section were produced with the help of the coding assistant Claude Code (Claude Opus 5.5, Anthropic); all results were checked by the authors.

\section{Outlook and related work}\label{sec:rel_work}

Our results assume that the structure graph $\mathcal{G}$ is completely known. Could this be relaxed to partial knowledge? For linear drifts and acyclic graphs, substantial information can be recovered without knowing the graph \cite{Amendola2025}. If generic identification suffices, i.e., identifiability for almost all values of a finite-dimensional parameter set, then steady-state data from one intervention per strongly connected component recovers Ornstein--Uhlenbeck parameters up to a global scaling factor, under conditions on the graph and the spectrum \cite{Salehkaleybar2026OneIntervention}. For nonlinear drifts, genericity is more subtle, but in principle experimental data can yield identifiability even without graph structure \cite{Inglese1997, botvinickgreenhouse2026uniquerecoverytransportmaps}.

An alternative source of information could be non-Gaussian driving noise, as demonstrated for linear drifts in \cite{recke2026identifiabilityestimationcontinuouslyapunov}. For nonlinear drifts, the Fokker--Planck equation \eqref{eq:FP} would be replaced by a nonlocal integro-differential equation. Finally, stationary distributions under certain parametric assumptions such as chemical reaction networks \cite{Grunberg2025} could reveal information on the graph, although this work, like ours, assumes a known network.

\ifprocl
\else
\section*{Acknowledgements}
We acknowledge support from the DAAD programme Konrad Zuse Schools of Excellence in Artificial Intelligence, sponsored by the Federal Ministry of Research, Technology and Space, from the European Research
Council (ERC) under the European Union’s Horizon 2020 research and innovation programme
(grant agreement No 883818), and from the German Federal Ministry of Education and Research and the Bavarian State Ministry for Science and the Arts. The authors take full responsibility for the contents of this article.
\fi

\bibliographystyle{amsalpha}
\bibliography{bibliography}

\newcommand{\etalchar}[1]{$^{#1}$}
\providecommand{\bysame}{\leavevmode\hbox to3em{\hrulefill}\thinspace}
\providecommand{\MR}{\relax\ifhmode\unskip\space\fi MR }
\providecommand{\MRhref}[2]{%
  \href{http://www.ams.org/mathscinet-getitem?mr=#1}{#2}
}
\providecommand{\href}[2]{#2}
\begin{thebibliography}{BGMY23}

\bibitem[ABHM25]{Amendola2025}
Carlos Améndola, Tobias Boege, Benjamin Hollering, and Pratik Misra,
  \emph{Structural identifiability of graphical continuous {L}yapunov models},
  2025, arXiv:2510.04985.

\bibitem[ASSP25]{aguirre2025gene}
Matthew Aguirre, Jeffrey~P. Spence, Guy Sella, and Jonathan~K. Pritchard,
  \emph{Gene regulatory network structure informs the distribution of
  perturbation effects}, PLoS Comput. Biol. \textbf{21} (2025), no.~9, Paper
  No. e1013387.

\bibitem[BGL14]{BGL2014Book}
Dominique Bakry, Ivan Gentil, and Michel Ledoux, \emph{Analysis and geometry of
  {M}arkov diffusion operators}, Grundlehren der mathematischen Wissenschaften,
  vol. 348, Springer, Cham, 2014.

\bibitem[BGMY23]{BotvinickGreenhouse2023}
Jonah Botvinick-Greenhouse, Robert Martin, and Yunan Yang, \emph{Learning
  dynamics on invariant measures using {PDE}-constrained optimization}, Chaos
  \textbf{33} (2023), no.~6, Paper No. 063152.

\bibitem[BGY26]{botvinickgreenhouse2026uniquerecoverytransportmaps}
Jonah Botvinick-Greenhouse and Yunan Yang, \emph{On the unique recovery of
  transport maps and vector fields from finite measure-valued data}, 2026,
  arXiv:2604.07671.

\bibitem[BKRS15]{Bogachev2015FPBook}
Vladimir~I. Bogachev, Nicolai~V. Krylov, Michael R\"ockner, and Stanislav~V.
  Shaposhnikov, \emph{Fokker-{P}lanck-{K}olmogorov equations}, Mathematical
  Surveys and Monographs, vol. 207, American Mathematical Society, Providence,
  RI, 2015.

\bibitem[BLD26]{bleile2026efficient}
Fabian Bleile, Sarah Lumpp, and Mathias Drton, \emph{Efficient learning of
  stationary diffusions with {S}tein-type discrepancies}, The 29th
  International Conference on Artificial Intelligence and Statistics, 2026.

\bibitem[BM26]{boeken2026causalgraphsmarkovproperties}
Philip Boeken and Joris~M. Mooij, \emph{Causal graphs, {M}arkov properties and
  do-calculus for stochastic differential equations}, 2026, arXiv:2607.12140.

\bibitem[Boc46]{Bochner1946}
S.~Bochner, \emph{Vector fields and {R}icci curvature}, Bull. Amer. Math. Soc.
  \textbf{52} (1946), 776--797.

\bibitem[DDK24]{Dettling2024Lasso}
Philipp Dettling, Mathias Drton, and Mladen Kolar, \emph{On the {L}asso for
  graphical continuous {L}yapunov models}, Proceedings of the Third Conference
  on Causal Learning and Reasoning, Proceedings of Machine Learning Research,
  vol. 236, PMLR, 2024, pp.~514--550.

\bibitem[DGFII25]{Gregorio2025DiffusionProcessRegression}
Alessandro De~Gregorio, Dario Frisardi, Stefano Iacus, and Francesco Iafrate,
  \emph{Adaptive elastic-net estimation for sparse diffusion processes}, Stat.
  Inference Stoch. Process. \textbf{28} (2025), no.~3, Paper No. 22.

\bibitem[DHA{\etalchar{+}}23]{Dettling2023}
Philipp Dettling, Roser Homs, Carlos Am\'endola, Mathias Drton, and
  Niels~Richard Hansen, \emph{Identifiability in continuous {L}yapunov models},
  SIAM J. Matrix Anal. Appl. \textbf{44} (2023), no.~4, 1799--1821.

\bibitem[EL00]{Elowitz2000}
Michael Elowitz and Stanislas Leibler, \emph{A synthetic oscillatory network of
  transcriptional regulators}, Nature \textbf{403} (2000), 335--338.

\bibitem[GDV23]{Grunberg2025}
Theodore~W. Grunberg and Domitilla Del~Vecchio, \emph{Identifiability of
  chemical reaction networks with intrinsic and extrinsic noise from stationary
  distributions}, SIAM J. Appl. Dyn. Syst. \textbf{22} (2023), no.~3,
  2206--2241.

\bibitem[Ing97]{Inglese1997}
Gabriele Inglese, \emph{A stochastic perturbation method for studying inverse
  problems of dynamics}, Appl. Anal. \textbf{64} (1997), no.~3-4, 391--398.

\bibitem[KM18]{Kolesnikov2018}
Alexander~V. Kolesnikov and Emanuel Milman, \emph{Poincar{\'e} and
  {B}runn-{M}inkowski inequalities on the boundary of weighted {R}iemannian
  manifolds}, Amer. J. Math. \textbf{140} (2018), no.~5, 1147--1185.

\bibitem[KS07]{Kar2007}
Sayan Kar and Soumitra Sengupta, \emph{The {R}aychaudhuri equations: A brief
  review}, Pramana \textbf{69} (2007), no.~1, 49--76.

\bibitem[Lim88]{Lima1988}
Elon~L. Lima, \emph{The {J}ordan-{B}rouwer separation theorem for smooth
  hypersurfaces}, Amer. Math. Monthly \textbf{95} (1988), no.~1, 39--42.

\bibitem[LKS24]{Lorch2024StationaryDiffusions}
Lars Lorch, Andreas Krause, and Bernhard Sch{\"{o}}lkopf, \emph{Causal modeling
  with stationary diffusions}, AISTATS 2024, Proceedings of Machine Learning
  Research, vol. 238, {PMLR}, 2024, pp.~1927--1935.

\bibitem[LT21]{Haesung2021FullSDEPackageUnderLocalIntegrability}
Haesung Lee and Gerald Trutnau, \emph{Existence, uniqueness and ergodic
  properties for time-homogeneous {I}t\^o-{SDE}s with locally integrable drifts
  and {S}obolev diffusion coefficients}, Tohoku Math. J. (2) \textbf{73}
  (2021), no.~2, 159--198.

\bibitem[LZ23]{Lacker2023StationaryDiffusionOnTree}
Daniel Lacker and Jiacheng Zhang, \emph{Stationary solutions and local
  equations for interacting diffusions on regular trees}, Electron. J. Probab.
  \textbf{28} (2023), Paper No. 4.

\bibitem[MPS90]{MalletParet1990}
John Mallet-Paret and Hal~L. Smith, \emph{The {P}oincar\'e-{B}endixson theorem
  for monotone cyclic feedback systems}, J. Dynam. Differential Equations
  \textbf{2} (1990), no.~4, 367--421.

\bibitem[RH26]{recke2026identifiabilityestimationcontinuouslyapunov}
Cecilie~Olesen Recke and Niels~Richard Hansen, \emph{Identifiability and
  estimation in continuous {L}yapunov models}, 2026, arXiv:2603.17142.

\bibitem[RS11]{Raulot2011}
S.~Raulot and A.~Savo, \emph{A {R}eilly formula and eigenvalue estimates for
  differential forms}, J. Geom. Anal. \textbf{21} (2011), no.~3, 620--640.

\bibitem[Sal26]{Salehkaleybar2026OneIntervention}
Saber Salehkaleybar, \emph{One intervention per component is enough: Towards
  identifiability in linear stochastic dynamics from steady state}, 2026,
  arXiv:2609.19955; presented at ICML 2026.

\bibitem[SBB25]{Shen2025}
Yunyi Shen, Renato Berlinghieri, and Tamara Broderick, \emph{Multi-marginal
  {S}chr{ö}dinger bridges with iterative reference refinement}, Proceedings of
  The 28th International Conference on Artificial Intelligence and Statistics,
  vol. 258, PMLR, 03--05 May 2025, pp.~3817--3825.

\bibitem[SD26]{schwank2026arxiv}
Richard Schwank and Mathias Drton, \emph{Non-parametric recovery of causal
  diffusion mechanisms from steady-state observations}, 2026, arXiv:2606.30467.

\bibitem[SSK{\etalchar{+}}21]{Song2021GenerativeSDE}
Yang Song, Jascha Sohl{-}Dickstein, Diederik~P. Kingma, Abhishek Kumar, Stefano
  Ermon, and Ben Poole, \emph{Score-based generative modeling through
  stochastic differential equations}, 9th International Conference on Learning
  Representations, {ICLR} 2021, OpenReview.net, 2021.

\bibitem[SW14]{Saumard2014LogConcaveReview}
Adrien Saumard and Jon~A. Wellner, \emph{Log-concavity and strong
  log-concavity: a review}, Stat. Surv. \textbf{8} (2014), 45--114.

\bibitem[VH20]{varando20a}
Gherardo Varando and Niels~Richard Hansen, \emph{Graphical continuous
  {L}yapunov models}, Proceedings of the 36th Conference on Uncertainty in
  Artificial Intelligence (UAI), PMLR, 03--06 Aug 2020, pp.~989--998.

\end{thebibliography}

\end{document}